\documentclass{article}
\usepackage{amssymb,amsmath,amsthm}
\usepackage{graphicx}
\usepackage[T1]{fontenc}
\usepackage{yfonts}

\newtheorem{theorem}{Theorem}

\newtheorem*{thmDH}{Theorem DH}

\newtheorem*{thmCoKer}{Theorem~\ref{thm:coker}}
\newtheorem*{thmTrans}{Theorem~\ref{thm:trans}}
\newtheorem*{thmSols}{Theorem~\ref{thm:smoothSols}}
\newtheorem*{thmGerms}{Theorem~\ref{thm:Germs}}
\newtheorem{lemma}{Lemma}
\newtheorem{proposition}{Proposition}

\newtheorem{corollary}{Corollary}
\newtheorem{definition}{Definition}

\def\bN{{\mathbb N}}

\def\bR{{\mathbb R}}

\def\bZ{{\mathbb Z}}

\def\cF{{\mathcal F}}
\def\cG{{\mathcal G}}

\def\cI{{\mathcal I}}

\def\cS{{\mathcal S}}
\def\cT{{\mathcal T}}

\def\fg{{\mathfrak g}}

\def\F{{_{F}}}
\def\Fn{{_{F_n}}}
\def\Fo{{_{F_1}}}
\def\G{{_{G}}}
\def\Gn{{_{G_n}}}

\def\FG{{_{FG}}}
\def\FnG{{_{F_nG}}}
\def\FoG{{_{F_1G}}}
\def\FoGo{{_{F_1G_1}}}
\def\FnGn{{_{F_n G_n}}}
\def\Rt{{\bR^2}}

\def\vf{{\mathfrak X}}
\def\loc{loc}

\def\id{\hbox{id}}

\def\gap{\mathop{gap}}

\def\coker{\hbox{\rm coker\,}}

\def\smooth{{C^\infty(\bR^2)}}

\def\smooth{C^\infty(\Rt)}
 
\title{Solvability of the cohomological equation \\for regular vector fields on the plane}

\author{Roberto De Leo}

\begin{document}
\maketitle

{\sl Math Subject Classification:} Primary: 37J99, 53C12, 35F05

{\sl Keywords:} Cohomological equation, foliations of the plane, Hamiltonian vector fields on the plane, linear first-order PDEs

{\sl Address:} Dipartimento di Matematica, Universit\`a di Cagliari, 09100 Cagliari, Italy

{\sl Email:} deleo@unica.it

%
\begin{abstract}
  We consider planar vector field without zeroes $\xi$ and study the image 
  of the associated Lie derivative operator $L_\xi$ acting on the space of 
  smooth functions.
  We show that the cokernel of $L_\xi$ is infinite-dimensional as soon as $\xi$ 
  is not topologically conjugate to a constant vector field and that, if the 
  topology of the integral trajectories of $\xi$ is ``simple enough'' 
  (e.g. if $\xi$ is polynomial) then $\xi$ is transversal to a Hamiltonian 
  foliation. We use this fact to find a large explicit subalgebra of the 
  image of $L_\xi$ and to build an embedding of $\Rt$ into $\bR^4$ which
  rectifies $\xi$. Finally we use this embedding to characterize the 
  functions in the image of $L_\xi$.
\end{abstract}
%
\section{Introduction, definitions and main results}
The study of planar vector fields has a long history going back to Poincar\`e and
Bendixson (see~\cite{God83} and~\cite{CC03} for more details and bibliography on 
this topic).
The topological classification of \emph{regular} (i.e. without zeros) vector fields 
on the plane was completed first by W.~Kaplan~\cite{Kap40,Kap48}, using an ad-hoc 
topological tool (chordal systems), based on previous works of his advisor 
H.~Whitney~\cite{Whi32,Whi33,Whi41}. In this paper we rather use the more general
concept of \emph{inseparable leaves} and \emph{separatrices}, introduced by 
L.~Markus~\cite{Mar54} while working at the extension of Kaplan's results to the 
more general problem of the topological classification of all planar vector fields.

\begin{figure}
  \centering
  \includegraphics[width=5cm]{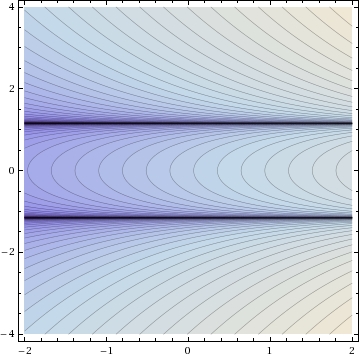}\hskip.75cm\includegraphics[width=5cm]{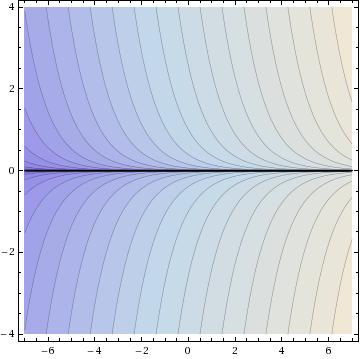}
  \caption{%
    \small
    Level sets of $F(x,y)=(y^2-1)e^x$ (left) and $G(x,y)=ye^x$ (right). 
    The first foliation has separatrices $y=\pm1$, the second has none.
  }
  \label{fig:s3}
\end{figure}

We recall a few standard basic concepts and definitions that will be used in the 
paper. We denote by $\vf_r(\Rt)$ the set of all smooth regular vector fields in the 
plane, by $\cF_\xi$ the foliation of the integral trajectories of $\xi\in\vf_r(\Rt)$
and by $\pi_\xi:\Rt\to\cF_\xi$ the canonical projection that sends every point 
in the leaf\footnote{Throughout the paper we refer to the points of $\cF_\xi$ as 
\emph{integral trajectories} or \emph{leaves} depending on the aspect of them we 
want to emphasize.}
passing through it. We endow $\cF_\xi$ with the canonical quotient topology.
It was shown by Haefliger and Reeb~\cite{HR57} that $\cF_\xi$ admits the
structure of a 1-dimensional simply connected second countable non (necessarily) 
Hausdorff smooth manifold; the smooth structure is characterized by the property 
that the restriction of $\pi_\xi$ to every transversal line $\ell$ is a diffeomorphism 
onto its image.
Two integral trajectories $s_i$, $i=1,2$, of $\xi$ are said \emph{inseparable} when their 
projections $\pi_\xi(s_i)$ cannot be separated in the topology of $\cF_\xi$ 
(e.g. see Fig.~\ref{fig:s3}). We denote by $\cI_{\cF_\xi,s}$ the set of all leaves 
distinct from $s$ inseparable from it (note that $\cI_{\cF_\xi,s}$ is empty for all but 
countably many leaves) and by $\cS_{\cF_\xi}$ the (countable) set of leaves for which 
$\cI_{\cF_\xi,s}$ is not empty.
A leaf $s$ is called a \emph{separatrix} when the boundary of every
neighbourhood of $\pi_\xi(s)$ contains more than two points. 
The set of all separatrices is the closure of $\cS_{\cF_\xi}$~\cite{Mar54}. 
In the present paper we will rather use the term separatrix to indicate just the
elements of $\cS_{\cF_\xi}$ since their limit points play no role in our work.
Every plane foliation is orientable and, correspondingly, to each set $\cI_{\cF_\xi,s}$ 
can be given a natural order; we say that two separatrices are \emph{adjacent} 
if they are next to each other with respect to this order.

We introduce now a few specific definition we will need throughout the paper.
\begin{definition}
  Two vector fields $\xi$ and $\xi'$ are \emph{strongly proportional} if they
  are proportional through a non-zero smooth function.
  A vector field $\xi$ is \emph{intrinsically Hamiltonian} if it is strongly
  proportional to a Hamiltonian vector field and is \emph{transversally Hamiltonian} 
  if it is transversal to a \emph{Hamiltonian foliation} $\cG$, i.e. to the level 
  sets of a regular smooth function $G$ (we say that $G$ is a Hamiltonian for $G$). 
\end{definition}
It is easily seen that a regular vector field is intrinsically Hamiltonian iff
the PDE $L_\xi f=0$ admits a regular smooth solution and is transversally 
Hamiltonian iff is is solvable the differential inequality $L_\xi f>0$.
\begin{definition}
  A foliation $\cF_\xi$ (or simply the vector field $\xi$) is of \emph{finite type} 
  if $\cS_{\cF_\xi}$ is closed and every set $\cI_{\cF_\xi,s}$ is finite. 
\end{definition}
In this case the complement of the set of separatrices 
is the disjoint union of countably many unbounded open sets named by 
Markus~\cite{Mar54} \emph{canonical regions} and the boundary of each canonical 
region has a finite number of connected components. We recall that examples
of smooth or even analytic foliations of the plane with a dense set of separatrices
are known in literature (see~\cite{Waz34} and~\cite{Mul76a}). While 
there are reasons to believe that such foliations are generic
in some ``combinatorial'' sense, the set of foliations of finite type is 
nevertheless of great importance since important natural categories of regular
vector fields leads to them. For example every polynomial vector field is 
of finite type: finite bounds for the number of the inseparable leaves 
of a polynomial vector field were find first by 
Markus~\cite{Mar72} and later improved independently by M.P.~Muller~\cite{Mul76b} 
and S.~Schecter and M.F. Singer~\cite{SS80}. It is easy to verify that are of
finite type also all regular vector fields strongly proportional to
those of the kind $\xi(x,y)=\left(a(y),b(y)\right)$, where $(a,b)$
is a generic pair of Morse functions of one variable (so that 
$a^2+b^2$ is strictly positive).
\begin{definition}
  A \emph{complete set of transversals} (CST) for $\cF_\xi$ is a set of 
  lines $\cT_\xi=\{\ell_i\}$, one for each separatrix of $\cF_\xi$, such that
  every $\ell_i$ is transversal to $\cF_\xi$ and cuts the corresponding 
  separatrix $s_i$ and the set $\{\pi_\xi(\ell_i)\}$ covers $\cF_\xi$.

  We call \emph{gap} of $g\in\smooth$ between two adjacent separatrices $s_1$ and $s_2$
  with respect to the CST $\cT_\xi$ the limit (if it exists) 
  $$
  \gap_{\cT_\xi}(g;s_1,s_2)=\lim_{p\to p_1}\int_0^{T_p} g\left(\Phi^t_\xi(p)\right)dt\,,   
  $$
  where the point $p\in\ell_1$ tends to $p_1=\ell_1\cap s_1$, $\Phi^t_\xi$ is the
  flux of $\xi$ and $T_p$ is the unique number s.t. 
  $\Phi^{T_p}_\xi(p)\in \ell_2$\footnote{Such number exists for $s_1$ and $s_2$ 
  are inseparable and is unique for every transversal cuts each leaf at most once.} 
\end{definition}
Finally we set a few notations on spaces of germs we are going to use in the 
last section. Let $a\in\bR$. We denote by $H^r_a$ the ring of left germs at $a$ 
of functions in $C^r(-\infty,a)$, i.e. the equivalence classes determined by 
the equivalence relation $h\simeq h'$ if $h$ and $h'$ coincide in some interval 
of the form $(a-\epsilon,a)$ for some $\epsilon>0$, and by $G^r_a$ the subring 
of the left germs in $H^r_a$ which can be extended to a continuous function at $a$
together with their derivatives up to order $r$.
Similarly, let $I=\{a\}\times[b_1,b_2]$ and set $L_I=(-\infty,a]\times\bR\setminus I$.
We denote by $H^r_I$ the ring of left germs at $I$ of functions of $C^r(L_I)$,
i.e. $h\simeq h'$ if $h$ and $h'$ coincide in some set $(U\cap L_I)\setminus I$, 
where $U$ is a neighbourhood of $I$, and by $G^r_I$ the subring of germs of functions 
of $H^r_I$ which can be extended to $C^r$ functions on the whole $L_I$.
\begin{definition}
  We call \emph{singular left germs at $a\in\bR$} the elements of the quotient 
  ring $SG^r_a=H^r_a/G^r_a$ and \emph{singular left germs at $I=\{a\}\times[b_1,b_2]$} 
  the elements of the quotient ring $SG^r_I=H^r_I/G^r_I$.
\end{definition}
%

Let us now turn to the topics of the present paper.
Recently S.P.~Novi\-kov~\cite{Nov08} (in case of smooth functions) and 
G.~Forni~\cite{For97} (for functional spaces of integrable functions) proved, as
a generalization of the well-known Diophantine phenomena in the torus, that the 
first order homogeneous partial differential operator associated to a generic 
vector field on a compact surface has an infinite-dimensional cokernel.
In Section~\ref{sec:coker} we generalize this result to the plane with the 
following result:
%
\begin{theorem}
  \label{thm:coker}
  Either $\xi$ is topologically conjugated to the constant vector field, 
  in which case $\dim\coker L_\xi=0$, or $\dim\coker L_\xi=\infty$.
\end{theorem}
%
%
It was proved by J.~Weiner~\cite{Wei88}, using a different terminology, 
that every Hamiltonian foliation is transversally Hamiltonian. 
In Section~\ref{sec:trans} we extend Weiner's result in the following way:
\begin{theorem}
  \label{thm:trans}
  Every plane foliation of finite type is transversally Hamiltonian.
\end{theorem}
%
%
In Section~\ref{sec:Im} we use Theorem~\ref{thm:trans} to
characterize the image of $L_\xi$ when $\xi$ is intrinsically Hamiltonian
or of finite type.
We recall that the question of the solvability of the so-called cohomological equation
\begin{equation}
  \label{eq:ce}
  L_\xi f=g\,,
\end{equation}
is of purely global nature: 
it is well known indeed that, if $\xi$ is regular, every point $p\in\Rt$ has a 
neighbourhood $U_p$ such that $L_\xi(C^\infty(U_p))=C^\infty(U_p)$.
A solution to~(\ref{eq:ce}) in $U_p$ is given explicitly by
$$
f(p) = \int_{0}^{T_p}g(\Phi^t_\xi(p_\ell))\; dt + h(\varphi(p))\,,
$$
where $\Phi^t_\xi$ is the flow of $\xi$, $T_p$ is the time needed to reach $p$ under
the action of $\Phi^t_\xi$ from the point $p_\ell$ lying on a fixed line $\ell$ embedded 
in $U_p$ and everywhere transversal to $\xi$, $\varphi\in C^\infty(U_p)$ is the 
(functional) generator of $\ker L_\xi$ and $h\in C^\infty(\bR)$.
The most general global result known for the action of a single regular vector field 
on the space of smooth functions is the following theorem by H\"ormander and 
Duistermaat~\cite{DH72}, which shows that a non-trivial behaviour of $L_\xi$ must correspond 
to a non-trivial structure of the foliation $\cF_\xi$ of the integral trajectories of $\xi$ 
and viceversa:
\begin{thmDH}
  Let $M$ be an open connected manifold and $\xi\in\vf_r(M)$.
  Then $L_\xi(C^\infty(M))=C^\infty(M)$ iff $\cF_\xi$ admits a global transversal\footnote{By global 
    transversal we mean a codimension-1 embedded submanifold of $M$ which is transversal 
    to $\cF_\xi$ at every point and cuts every leaf exactly once.}.
\end{thmDH}
For the case of regular vector fields in $\Rt$ we show (see Proposition~\ref{thm:xi'} and
the paragraph below it) that for every $\xi$ intrinsically 
Hamiltonian or of finite type there exist two commuting vector fields $\xi'_\F$ and $\xi'_\G$,
the first strongly proportional to $\xi$ and the second transversal to it ($\xi'_\G$ diverges
on some separatrices in the finite type case but this does not hinder the result), for which 
the following holds:
\begin{theorem}
  \label{thm:smoothSols}
  A function $g\in\smooth$ belongs to $L_{\xi'_\F}(\smooth)$ iff all functions 
  $L^k_{\xi'_\G} g$, $k\in\bN$, have finite gap for all pairs of adjacent separatrices 
  of $\xi$.
\end{theorem}
Finally, in Section~\ref{sec:graph} we show that there exists an embedding $\hat\Phi_\FG$ 
of $\Rt$ into $\bR^4$ which rectifies globally both $\xi'_\F$ and $\xi'_\G$ at the same time. 
This setting for the cohomological equation is dual to the initial one in the following 
sense: in the original one the ambient space is always the same and the qualitative 
difference between inequivalent cohomological equations resides in the analytic 
expression of the vector field; in the embedding, instead, the cohomological equation 
has always the same analytic expression and is the geometry of the ambient 
space $\hat\Phi_\FG(\Rt)$ that determines the solvability of the equation.

In this second setting we prove the following results:
\begin{theorem}
  \label{thm:Germs}
  There exists a countable family of intervals $I_j=\{a_j\}\times[b_{j,1},b_{j,2}]$ 
  and of ring homomorphisms 
  $\theta^{(r)}_j:SG^r_{I_j}\to SG^r_{a_j}$ such that $g\in L_{\xi'_\F}(C^r(\Rt))$ iff 
  $[(\hat\Phi_\FG)_* g]_{SG_{I_j}}\in\ker\theta^{(r)}_j$ for all $\theta^{(r)}_j$.
\end{theorem}
%
%

Note that in this paper we are interested only to the action of $\xi$ on smooth functions; 
concerning the global solvability in other functional spaces, e.g. of entire functions or 
Gevrey-type functions in the realm of global Cauchy-Kowalevs\-ka\-ya theorem see~\cite{GM03,GG07} 
and the references therein. Note also that there is a qualitative difference between 
the case of a single operator $L_\xi$ and the case of two or more operators 
$\{L_{\xi_1},\dots,L_{\xi_k}\}$; it has been shown indeed by M.~Gromov~\cite{Gro86} that, 
on every smooth manifold $M$,
$$
L_{\xi_1}(C^\infty(M))+L_{\xi_2}(C^\infty(M)) = C^\infty(M)
$$
for any pair of vector fields in mutual generic position.
%
%
\section{$\coker L_\xi$}
\label{sec:coker}
As pointed out above in Theorem~DH, if $\cF_\xi$ admits a global transversal 
the method of characteristics provides a global solution to the cohomological
equation (\ref{eq:ce}) for every $g\in\smooth$, so that $L_\xi(\smooth)=\smooth$ 
and $\coker L_\xi=\{0\}$.
The obstruction to the existence of global transversals is the presence of 
separatrices since no smooth line $\ell$ can, at the same time, be transversal 
to $\cF_\xi$ and intersect any pair leaves inseparable from each other. 

In absence of global transversals, one can try to solve $L_\xi f=g$ recursively 
in the following way. Let $s$ be a separatrix for $\xi$ and denote by $\ell$ any 
transversal through it and by $U_\ell=\pi^{-1}_\xi(\ell)\subset\Rt$ the saturated open 
set containing $\ell$. Since $U_\ell$ is a proper subset of $\Rt$, its boundary 
is non-empty and equal to the union of the sets $\cI_{\cF_\xi,\tilde s}$ corresponding to
all leaves $\tilde s$ cut by $\ell$.
By construction $\xi$, once restricted to $U_\ell$, admits a global transversal 
(the line $\ell$) and therefore $L_\xi(C^\infty(U_\ell))=C^\infty(U_\ell)$. 
Let now $g_\ell$ be any solution, in $U_\ell$, to $L_\xi f=g$. We can try to extend 
$g_\ell$ beyond $U_\ell$ by selecting any boundary component $s'$ of $\partial U_\ell$ 
and any transversal $\ell'$ passing through it. The function $g_\ell$ restricts 
on $\ell'\cap U_\ell$ to a smooth function $\hat g_{\ell'}$; if we can 
extend $\hat g_{\ell'}$ to a smooth function $g_{\ell'}$ defined on the whole $\ell'$ 
then, via the method of characteristics applied to the set $U_{\ell'}=\pi^{-1}_\xi(\ell')$ 
and using $g_{\ell'}$ as initial condition on $\ell'$, we can smoothly extend $g_\ell$ 
to $U_{\ell'}$. Assuming that one can always extend a local solution across 
transversals as described above, proceeding recursively until no 
separatrices are left we end up with a global solution to (\ref{eq:ce}).

We are going to use the gap to provide a quantitative criterion for the
existence of continuous solutions. 
While the gap of a function clearly depends on the CST chosen, whether 
it exists and is bounded does not:
\begin{proposition}
  If the gap of $g\in\smooth$ between two adjacent separatrices $s_1$ and $s_2$ with 
  respect to a CST $\cT_\xi$ exists and it is finite, then it exists and it 
  is finite also with respect to every other CST $\cT'_\xi$.
\end{proposition}
\begin{proof}
  Let $\ell'_1,\ell'_2\in\cT'_\xi$ be the two transversal to $s_1$ and $s_2$ in the
  second CST. Then
  $$
  \gap_{\cT'_\xi}(g;s_1,s_2) = \gap_{\cT_\xi}(g;s_1,s_2) 
  + A_1 + A_2
  $$
  for
  $$
  A_1 = \int_{p'_1}^{p_1} g\left(\Phi^t_\xi(p'_1)\right)dt\,,\; 
  A_2 = \int_{p_2}^{p'_2} g\left(\Phi^t_\xi(p_2)\right)dt
  $$
  where the integral defining $A_i$, $i=1,2$, is evaluated along $s_i$.
  Recall that, due to the method of characteristics, the values on a leaf of a 
  local solution to the cohomological equation are completely determined by 
  the value of the solution in any point of the leaf and they are finite on
  the whole leaf iff they are finite at a single point. Hence, if the gap of
  $g$ between $s_1$ and $s_2$ with respect to $\cT_\xi$, both $A_i$ are finite
  since they are given by integrals of bounded functions over compact sets. 
\end{proof}
It is already implicit in the previous proof that the existence and boundedness 
of the gap of a function $g$ is related to the extendability of local solutions
of the cohomological equation having $g$ as rhs. Below we prove this fact and
then use it to prove the main result of the section.
%
%
\begin{proposition}
  \label{thm:gCont}
  A global continuous solution to $L_\xi f=g$ exists iff $g$ has finite gap between 
  every pair of adjacent separatrices of $\cF_\xi$.
\end{proposition}
\begin{proof}
  We point out first that a continuous solution to $L_\xi f=g$, $g\in\smooth$, is much more 
  regular than it sounds since all such solutions are, by definition, smooth in the 
  $\xi$ direction. In particular the integral of $df$ along the integral trajectories
  of $\xi$ is well-defined even for continuous solutions of~(\ref{eq:ce}) since the
  restriction of $df$ on these integral trajectories depends only on $L_\xi f$.


  The condition in the hypothesis of the theorem is clearly necessary for, if a continuous 
  solution $f$ exists, then for a given $\cT_\xi$ we have
  $$
  \gap_{\cT_\xi}(g;s_1,s_2) = \lim_{p\to p_1}\int_0^{T_p} df = f(p_2)-f(p_1)\,.
  $$
  Note that the gap of $g$ between $s_1$ and $s_2$ depends only on the intersection 
  of the two separatrices with the relative transversals in $\cT_\xi$.

  Now assume that a solution $f_1$ is defined in $U_1=\pi^{-1}(\ell_1)$ and that the gap 
  of $g$ between $s_1$ and $s_2$ is finite. Then the restriction of $f_1$ on $\ell_2$ can
  be continued to a continuous function on the whole $\ell_2$ and therefore, via the
  the method of characteristics, to the whole $U_2=\pi^{-1}(\ell_2)$. 
  The new function $f_2$ defined on $U_1\cup U_2$ coincides, by construction, with
  $f_1$ in $U_1\cap U_2$, is continuous in $U_1\cup U_2$
  and clearly does not dependent on the choice of the particular CST used in the extension.
  By proceeding recursively until all separatrices are taken into account we 
  end up with a global continuous solution to~(\ref{eq:ce}).
\end{proof}
We are now in condition to prove Theorem~\ref{thm:coker}.
%
\begin{thmCoKer}
  If $\xi$ has at least a pair of separatrices then $\dim\coker L_\xi=\infty$.
\end{thmCoKer}
%
\begin{proof}
  We can assume without loss of generality that $\xi$ is complete\footnote{This 
  is true for any smooth vector field on a manifold, e.g. see \cite{God83}, Proposition~1.19;
  in this case, since $\xi$ is regular, we could simply assume that it has
  unitary Euclidean length.}.
  Under this assumption the gap of every non-zero constant function is 
  infinite for it is proportional to $T_p$,
  which clearly diverges for $p\to p_1$. Then the gap diverges also on every 
  function which is minored by a non-zero constant, e.g. the polynomials 
  $p_{n,m}(x,y)=1+x^{2n}+y^{2m}$, so that the image of $L_\xi$ misses infinitely 
  many linearly independent functions, i.e. $\dim\coker L_\xi=\infty$. 
\end{proof}
%

%
\section{$L_\xi f>0$}
\label{sec:trans}
Finding criteria to characterize functions belonging to the image of $L_\xi$
is hard and in the case of a generic regular vector field we cannot state much more
than the fact that a necessary condition (but far from being sufficient) to belong 
to it is to have finite gap between all pairs of adjacent separatrices.
%
%
More can be said for the vector fields which are transversally Hamiltonian, 
which makes crucial studying the solvability of the differential inequality $L_\xi f>0$.
%
\begin{proposition}
  \label{thm:transvHam}
  Let $\xi\in\vf_r(\Rt)$, $\Omega_0=dx\wedge dy$ and $\omega_\xi=i_\xi\Omega_0$. 
  The following conditions are equivalent:
  \item 1. $\cF_\xi$ is transversally Hamiltonian;
  \item 2. the inequality $L_\xi f>0$ has a smooth solution;
  \item 3. $\omega_\xi\wedge df$ is a volume form for some $f\in\smooth$.
\end{proposition}
\begin{proof}
  Let $\cG$ be a Hamiltonian foliation transversal to $\cF_\xi$ and $G$ a Hamiltonian
  for $\cG$. Since $T\cG=\ker dG$ we must have $dG(\xi)\neq0$ at every point, so
  that either $L_\xi G>0$ or $L_\xi(-G)>0$ and viceversa. 
  Part 3 is due to the fact that $\omega_\xi\wedge dG=i_\xi dG\,\Omega_0=L_\xi G\,\Omega_0$.
\end{proof}
As mentioned in the introduction, Weiner~\cite{Wei88} proved that every intrinsically 
Hamiltonian falls in this class. Below, after a preparatory Lemma, we extend this result 
to every $\xi$ of finite type. 
%
\begin{lemma}
  Let $\xi$ be a regular vector field of finite type.
  Then $\cF_\xi$ admits a CST with the following
  property: for each separatrix $s\in\cS$, the saturated open set $\pi^{-1}_\xi(\pi_\xi(\ell))$ 
  of all leaves cutting the corresponding transversal $\ell\in\cT$ is equal to the union 
  of $s$ with the two canonical regions having $s$ as boundary component.
\end{lemma}
\begin{proof}
  Let $s$ be a separatrix, $U$ one of the two canonical regions having $s$ as boundary,
  $\ell$ the corresponding transversal in $\cT$ and $\ell_U$ the connected component of
  $\ell\setminus s$ which intersects $U$.
  Since $U$ admits a global transversal, there is a natural diffeomorphism $\psi$ of $U$ 
  into $\bR$ sending the leaves of $\cF_\xi$ into vertical lines. 
  If $\pi^{-1}_\xi(\pi_\xi(\ell_U))\neq U$ there is no geometrical obstruction to make 
  $\psi(\ell_U)$ either shorter or longer in the horizontal direction while keeping it 
  transversal to the vertical direction and without modifying it close to $s$ so that the first 
  projection of $\psi(\ell)$ on the first factor is surjective. After we do the same
  on the second canonical region $V$ we are left with a new transversal $\ell'$ such that
  $\pi^{-1}_\xi(\pi_\xi(\ell'))= U\cup V\cup s$.
\end{proof}
\begin{thmTrans}
  \label{thm:Xf>0}
  Every regular vector field of finite type is transversally Hamiltonian.
\end{thmTrans}
\begin{proof}
  We can assume without loss of generality that $\xi$ is complete and denote 
  by $\cT_\xi$ any CST having the property described in the Lemma above.
  The collection of open subsets $V_{s,i}$ defined by
  $$
  V_{s,i}=\{\Phi_\xi^t(\ell_s)\,|\,t\in(i,i+1)\}\,,\;s\in\cS_\xi\,,\;i\in\bZ\,,
  $$
  where $\Phi_\xi$ is the flow of $\xi$ and $\ell_s$ the transversal associated to $s$
  in $\cT_\xi$, is a locally finite open cover of $\Rt$. Indeed by hypothesis the union 
  of the $\pi_\xi(\ell_i)$ covers $\cF_\xi$ and therefore under the flow $\Phi_\xi$ the 
  $\ell_i$ cover the whole plane. Moreover, since the boundary of every canonical region 
  has only finitely many components, only finitely many of the $V_{s,i}$ cover any given
  point.

  Inside each $V_{s,i}$ every point $p$ can be written as $\Phi_\xi^t(q)$ for some 
  $q\in \ell_s$ so that we can define a smooth function $f_{s,i}(\Phi_\xi^t(q))=\phi(t)$, 
  where $\phi$ is any smooth function strictly monotonic for $t\in(0,1)$ and such that 
  $\phi|_{(-\infty,0)}\equiv0$ and $\phi|_{(1,\infty)}\equiv1$.
  Since each $V_{s,i}$ divides the plane in two connected components, 
  each $f_{s,i}$ can be extended to a smooth function on the whole plane by setting it identically
  to 1 in the component containing $\Phi_\xi^1(\ell_s)$ and identically 0 in the other.
  A direct calculation shows that $L_\xi f_{s,i}(p)=\phi'(t)>0$ within each $V_{s,i}$ while
  $L_\xi f_{s,i}$ is identically zero outside of it. Now recall that the set $\cS_\xi\times\bZ$ 
  is countable and let $n_{s,i}$ be any bijection of it with $\bN$. The series 
  $$
  f=\sum_{s\in\cS_\xi,i\in\bZ} 2^{-n_{s,i}} f_{s,i}
  $$
  converges to a continuous function (because the $f_{s,i}$ are uniformly bounded) which 
  is actually smooth because the derivatives of all positive orders of the $f_{s,i}$ have 
  compact support. By construction $L_\xi f\geq0$ but the inequality is strict because
  for every $x_0$ there exists at least one index $(s_0,i_0)$ such that $L_\xi f_{s_0,i_0}>0$.
\end{proof}
%
Note that the inequality $L_\xi f>\epsilon$, with $\epsilon>0$, requires stricter 
conditions to be solvable no matter how small $\epsilon$ is. E.g. it admits no 
smooth solutions if $\xi$ is complete for in that case, as pointed out in the 
previous section, all gaps of the constant function $\epsilon$ (and, a fortiori, 
all gaps of every function not smaller than it) would be infinite.
%
\section{$L_\xi(\smooth)$}
\label{sec:Im}
%
%
%
From this point on we will assume that $\xi$ is transversally Hamiltonian 
and we will denote by $F\in\smooth$ a generator of $\ker L_\xi$, 
so that $\ker L_\xi=F^*\left(C^\infty(\bR)\right)$, by $\cG$ the Hamiltonian foliation
transversal to $\cF_\xi$ and by $G$ any Hamiltonian of $\cG$.

A fundamental tool in our analysis will be the map $\Phi_\FG:\Rt\to\Rt$ 
defined by $x'=F(x,y)$, $y'=G(x,y)$. 
Assume first that $\xi$ is intrinsically Hamiltonian, so that $F$ is regular. 
In this case $\Phi_\FG$ is an immersion, since also $G$ is regular and the 
level sets of $F$ and $G$ are everywhere transversal by hypothesis, so 
that it induces on the source space the following metric and symplectic structures:
$$
g_\FG=\Phi_\FG^*((dx')^2+(dy')^2)=(dF)^2+(dG)^2\,,\;\Omega_\FG=\Phi_\FG^*(dx'\wedge dy')=dF\wedge dG\,.
$$
In particular $\Phi_\FG$ induces on the source space complex structure $J_\FG$,
whose real and imaginary spaces are $T\cF_\xi$ and $T\cG$, and a Poisson 
structure $\{,\}_\FG$.
Via $\Omega_\FG$ we can build a pair of commuting regular vector fields respectively 
tangent to $\cF_\xi$ and $\cG$. Recall that the Hamiltonian vector field $\eta$ associated
to a Hamiltonian $H$ with respect to a symplectic form $\Omega$ is defined by the
relation $i_{\eta}\Omega=dH$; in this case we write, with a slight abuse of notation,
that $\eta=\Omega^{-1}(dH)$.
\begin{proposition}
  \label{thm:xi'}
  Let $\xi'_\F=-\Omega_\FG^{-1}(dF)$, $\xi_\F=-\Omega_0^{-1}(dF)$, $\xi'_\G=\Omega^{-1}_\FG(dG)$ 
  and $\xi_\G=\Omega^{-1}_0(dG)$. The following relations hold:
  \begin{enumerate}
  \item $\Omega_\FG=(L_{\xi_\F}\!\!G)\; \Omega_0$.
  \item $\xi'_\F=\frac{1}{L_{\xi_\F}\!\!G}\xi_\F$, 
        $\xi'_\G=\frac{1}{L_{\xi_\F}\!\!G}\xi_\G$.
  \item $L_{\xi'_\F}F=0\,,\;L_{\xi'_\F}G=1\,,\;L_{\xi'_\G}F=1\,,\;L_{\xi'_\G}G=0\,.$
  \item $(\Phi_\FG)_*(\xi'_\F)=\partial_{y'}$ and $(\Phi_\FG)_*(\xi'_\G)=\partial_{x'}$ within $\Phi_\FG(\Rt)$.
  \item $\{F,G\}_{\FG}=L_{\xi_\F}G=1$.
  \item $[\xi'_\F,\xi'_\G]=0$.
  \item The pair $(\xi'_\F,\xi'_\G)$ is an orthonormal base for $g_\FG$.
  \item $L_{\xi'} g_\FG = L_{\eta'} g_\FG = 0$.
  \item $L_{\xi'}\Omega_\FG=L_{\eta'}\Omega_\FG=0$.
  \item $J_\FG\xi'_\F=\xi'_\G\,,\,\,J_\FG\xi'_\G=-\xi'_\F$.
  \end{enumerate}
\end{proposition}
\begin{proof}
  \begin{enumerate}
  \item A direct calculation show that 
    $\xi_\F=-\partial_y F\partial_x+\partial_x F\partial_y$, so that  
    $dF\wedge dG=(\partial_x F\partial_y G-\partial_y F\partial_x G)dx\wedge dy=(L_{\xi_\F}G)\Omega_0$.
  \item It is a direct consequence of the definition of $\xi'_\F$ and $\xi'_\G$ and (1).
  \item It is a direct consequence of $(2)$.
  \item Since $\Phi_\FG$ is not an injection, the push-forward of a vector field 
    $(\Phi_\FG)_*(\zeta)=T\Phi_\FG\circ\zeta\circ\Phi_\FG^{-1}$ is not well-defined unless 
    $T\Phi_\FG(\zeta)$ takes the same value on all points of $\Phi_\FG^{-1}(p)$ for 
    every $p\in\Phi_\FG(\Rt)$. This is the case for $\xi'_\F$ and $\xi'_\G$ since
    we get in both cases a constant vector field:
    $$((\Phi_\FG)_*(\xi'_\F))(x')=\xi'_\F(\Phi_\FG^*(x'))=\xi'_\F(F)=0$$ 
    $$((\Phi_\FG)_*(\xi'_\F))(y')=\xi'_\F(\Phi_\FG^*(y'))=\xi'_\F(G)=1$$ 
    and similarly for $\xi'_\G$.
  \item $\{F,G\}_\FG=\{\Phi_\FG^*x',\Phi_\FG^*y'\}_\FG=\Phi_\FG^*\{x',y'\}_0=\Phi_\FG^*1=1$.
  \item $[\xi'_\F,\xi'_\G]=[-\Omega_\FG^{-1}(dF),\Omega_\FG^{-1}(dG)]=\Omega_\FG^{-1}(\{F,G\}_\FG)=\Omega_\FG^{-1}(1)=0$.
  \item $g_\FG(\xi'_\F,\xi'_\F) = \left(dF(\xi'_\F)\right)^2+\left(dG(\xi'_\F)\right)^2=(L_{\xi'_\F}F)^2+(L_{\xi'_\F}G)^2=0+1$
  and similarly for the other combinations.
  \item It is a direct consequence of the previous item.
  \item This just restates that $\xi'_\F$ and $\xi'_\G$ are Hamiltonian with respect to $\Omega_\FG$.
  \item It is due to the fact that both $g_\FG$ and $\Omega_\FG$ are in canonical form with respect to $\xi'_\F$ and $\xi'_\G$.
  \end{enumerate}
\end{proof}
%
%
%
When $\xi$ is not intrinsically Hamiltonian but is of finite type then $F$ is not 
globally regular but nevertheless its differential goes to zero only on some of the separatrices,
so that the restriction of $\Phi_\FG$ to each of the canonical regions of $\xi$ is still an immersion.
Correspondingly, the pair of commuting regular vector fields $\xi'_\F$ and $\xi'_\G$ is well
defined within the canonical regions but, while $\xi'_\F$ is globally well-defined, $\xi'_\G$ 
diverges on the separatrices where $dF$ is zero. Note that there is no way to find a 
global substitute for $\xi'_\G$:
\begin{proposition}
  Let $\cF$ be a plane foliation of finite type. Then a pair of commuting regular 
  linearly independent vector fields $(\xi,\eta)$, with $\xi$ tangent to $\cF$, 
  exists iff $\cF$ is Hamiltonian.
\end{proposition}
\begin{proof}
  We showed in previous proposition that such pair always exists if $\cF$ is Hamiltonian.
  Assume then that it is not. In this case we can always find a smooth function $F$ with no maxima
  or minima whose differential vanishes on some of the separatrices and whose level sets 
  are the leaves of $\cF$ and a second function $G$, this one regular on the whole plane, 
  whose level sets are always transversal to $\cF$. Correspondingly we can always find
  two vector fields $\xi$ and $\eta$ s.t. 
  $$
  L_\xi F=0\,,\;L_\xi G=1\,,\;L_\eta G=0\,,\;L_\eta F\geq0\,.
  $$
  Let now $\alpha$ e $\beta$ the two smooth functions s.t. $[\xi,\eta]=\alpha\xi+\beta\eta$.
  Then 
  $$
  \alpha=\alpha L_\xi G+\beta L_\eta G = L_{[\xi,\eta]}G = L_\xi(L_\eta G)-L_\eta(L_\xi G)=0
  $$
  while 
  $$
  \beta L_\eta F=\alpha L_\xi F + \beta L_\eta F = L_{[\xi,\eta]}F = L_\xi(L_\eta F)-L_\eta(L_\xi F)=L_\xi(L_\eta F)
  $$
  namely $\beta=L_\xi[\log L_\eta F]$. Since $[\xi,\eta]$ has only the $\eta$ component, 
  the only thing we can do to make the commutator vanish is multiplying $\eta$ by some 
  non-zero factor $\mu$ since any other change would just introduce a $\xi$ component.
  On the other side
  $$
  [\xi,\mu\,\eta]=L_\xi\mu\,\eta+\mu[\xi,\eta]=L_\xi\mu\,\eta+\mu\,\beta\,\eta
  $$
  leading to $\mu=1/L_\eta F$; this function though is not smooth because the 
  differential of $F$ vanishes on some of the separatrices.
\end{proof}
Let us turn now to the study of the image of $L_\xi$. 
This is clearly equivalent to studying the image of $L_{\xi'_\F}$ but the latter 
is more convenient for the following two propositions:
\begin{proposition}
  \label{thm:Fy=G}
  The cohomological equation $L_{\xi'_\F}f(x,y)=g(x,y)$, restricted to the subalgebra
  $\Phi_\FG^*\left(\smooth\right)=\{\Phi_\FG^*\hat f\,|\,\hat f\in\smooth\}$, 
  writes, in the image of $\Phi_\FG$, as  
  \begin{equation}
    \label{eq:Fy=G}
    \frac{\partial\phantom{y'}}{\partial y'}\hat f(x',y')=\hat g(x',y')
  \end{equation}
  where $\hat f=(\Phi_\FG)_*f$ and $\hat g=(\Phi_\FG)_*g$.
\end{proposition}
\begin{proof}
  In general $\Phi_\FG$ is not injective so that, while the pull-back of function 
  $\Phi_\FG^*\hat f:=\hat f\circ\Phi_\FG$ is well-defined on the whole $\smooth$, the push
  forward $(\Phi_\FG)_*f:=f\circ\Phi_\FG^{-1}$ leads to a well-defined function only 
  within the subalgebra $\Phi_\FG^*\left(\smooth\right)$.
  Then from point (3) of Proposition~\ref{thm:xi'} follows that
  $$
  (\Phi_\FG)_*\left(L_{\xi'_\F}\left(\Phi_\FG^*\hat f\right)\right)
  = L_{(\Phi_\FG)_*{\xi'_\F}}\left((\Phi_\FG)_*\Phi_\FG^*\hat f\right)
  = \frac{\partial\phantom{y'}}{\partial y'}\hat f
  $$
\end{proof}
\begin{thmSols}
  A function $g\in\smooth$ belongs to $L_{\xi'_\F}(\smooth)$ iff all functions 
  $L^k_{\xi'_\G} g$, $k\in\bN$, have finite gap for all pairs of adjacent separatrices 
  of $\xi'_\F$.
\end{thmSols}
\begin{proof}
  As we already pointed out, every continuous solution to $L_{\xi'_\F} f=g$ is automatically
  smooth in the $\xi'_\F$ direction, i.e. $L^k_{\xi'_\F} f$ is continuous for every $k\in\bN$.

  Assume first that $\xi'_\F$ is intrinsically Hamiltonian.
  Since $\xi'_\F$ and $\xi'_\G$ commute and are globally well-defined, the first derivative 
  in the $\xi'_\G$ direction satisfies the cohomological equation 
  $L_{\xi'_\F}\left(L_{\xi'_\G} f\right)=L_{\xi'_\G} g$ and analogously the $k$-th derivative in 
  the ${\xi'_\G}$ direction satisfies $L_{\xi'_\F}\left(L^k_{\xi'_\G} f\right)=L^k_{\xi'_\G} g$. 
  Now we can use the claim of Lemma~\ref{thm:gCont} to conclude that each $L^k_{\xi'_\G} f$ is 
  globally continuous iff $L^k_{\xi'_\G} g$ has finite gap between every pair of adjacent
  separatrices.

  Assume now that $\xi'_\F$ is of finite type, so that ${\xi'_\G}$ is only well-defined within 
  the canonical regions of $\xi'_\F$. By repeating the same kind of arguments used in 
  Lemma~\ref{thm:gCont} it is clear that we can extend a smooth solution within a 
  saturated open set to the whole plane iff the gap of $L^k_{\xi'_\G} g$ has finite gap 
  between every pair of adjacent separatrices. Note indeed that in the definition of gap
  the values of $\xi'_\G$ on the separatrices are never used so the fact that $\xi'_\G$ 
  diverges on some of them does not hinder the evaluation of the gap.
\end{proof}
From Proposition~\ref{thm:Fy=G} and the surjectivity of $\partial_{y'}$ we 
get a large explicit subalgebra of the image of $L_{\xi'_\F}$:
\begin{proposition}
  \label{thm:subset}
  $\Phi_\FG^*\left(\smooth\right) \subset L_{\xi'_\F}(\smooth)$
\end{proposition}
This fact corresponds to two elementary observations: one, algebraic, that
$$
L_{\xi'_\F}\hat f(F,G)=L_{\xi'_\F}F\,\partial_{x'} \hat f(F,G) + L_{\xi'_\F}G\,\partial_{y'} \hat f(F,G) 
= \partial_{y'} \hat f(F,G)\,;
$$
the other, geometric, that the constant vertical vector field $\partial_{y'}$ 
on $\Phi_\FG(\Rt)$ can always be extended to the whole plane, where it is
surjective on $\smooth$.
%
\section{Local behaviour of functions of $L_\xi(\smooth)$ close to a pair 
of adjacent separatrices}
\label{sec:graph}
Proposition~\ref{thm:Fy=G} shows that locally, in the image of the map $\Phi_\FG$, 
the cohomological equations relative to vector fields $\xi'_\F$ look all the same, 
independently on the topology of their leaf spaces; 
the qualitative difference between them resides rather in the global geometry 
of the map $\Phi_\FG$. It is easy to verify that, as soon as $\xi'_\F$ has at 
least two pairs of separatrices, $\Phi_\FG$ cannot be injective, which is not 
optimal for several reasons. We bypass this problem by considering the map 
$\hat\Phi_\FG:\Rt\to\bR^4$ defined by $\hat\Phi_\FG(x,y)=(x,y,F(x,y),G(x,y))$. 
By construction $\hat\Phi_\FG$ is a diffeomorphism between $\Rt$ and 
$\Gamma_\FG=\hat\Phi_\FG(\Rt)\subset\bR^4$, the graph of $\Phi_\FG$. 
The symplectic, metric and almost complex structures
determined on $\Rt$ by $F$ and $G$, as pointed out at the beginning of the 
previous section, induce the same structures on $\Gamma_\FG$ via the
push-forward $(\hat\Phi_\FG)_*$.
We use on $\bR^4=\Rt\times\Rt$ coordinates $(x,y,x',y')$
and denote by $\pi_1$ and $\pi_2$ the projections on the first and second
factor. By definition $\pi_1\circ\hat\Phi_\FG=\id_\Rt$ and 
$\pi_2\circ\hat\Phi_\FG=\Phi_\FG$, so $\Gamma_\FG$ admits $(x,y)$ as global 
coordinates and $(F,G)$ as local coordinates at every point.
A direct calculation shows that 
$$
(\hat\Phi_\FG)_*(\xi'_\F)=\xi'_\F\oplus\partial_{y'}\,,\;
(\hat\Phi_\FG)_*(\xi'_\G)=\xi'_\G\oplus\partial_{x'}\,.
$$
In particular the projection on the second factor of the image of 
the leaves of $\cF_\xi$ and $\cG$ are, respectively, vertical and horizontal 
straight lines in the plane $(x',y')$.
All leaves which are inseparable one from the other are mapped to
disjoint open intervals of the same line, so that the images 
in the graph of any pair of adjacent separatrices of $\xi'_\F$ are 
separated by a vertical closed bounded interval $I$.
\begin{proposition}
  \label{thm:nc}
  For every pair of separatrices $s_1$ and $s_2$ of $\xi'_\F$, with
  $a=F|_{s_1\cup s_2}$, there exists a saturated open neighbourhood $U$ of 
  $s_1$ and $s_2$ on which $\Phi_\FG$ is injective and 
  $\Phi_\FG(U\cap \Phi_\FG^{-1}((a_1,a_2)\times(c_1,c_2)))=(a_1,a_2)\times(c_1,c_2)\setminus R$, 
  where $R=[a,a_2)\times[b_1,b_2]$ or $R=(a_1,a]\times[b_1,b_2]$, both $a_i$ 
  and $c_i$ can be infinite and $c_1<b_1\leq b_2<c_2$.
\end{proposition}
\begin{proof}
  Let $p_i\in s_i$, $i=1,2$, be any two points on the two separatrices,
  set $c_i=G(p_i)$ and denote by $\ell_i$ be the two leaves of $\cG$ 
  passing through the $p_i$. The two numbers $c_1$ and $c_2$ cannot be equal 
  since the restriction of $G$ to any leaf of $\cF_{\xi'_\F}$ is strictly 
  monotonic and, because of the inseparability of $s_1$ and $s_2$, there 
  are leaves of $\cF_{\xi'_\F}$ 
  cutting both $\ell_1$ and $\ell_2$; in particular $G(s_1)\cap G(s_2)=\emptyset$.
  Assume that $c_1<c_2$ (otherwise switch
  the names of the points), set $U_i=\pi_\xi^{-1}(\ell_i)$, $i=1,2$ and denote 
  by $V$ and $\Lambda$ respectively the union and intersection of $U_1$ and $U_2$.

  Assume first that $\Lambda$ is contained in $F<a$.
  We claim that the restriction of $\Phi_\FG$ to $V$ is injective. 
  Indeed let $A_i=U_i\setminus\Lambda$, $i=1,2$, so that $V=\Lambda\sqcup A_1\sqcup A_2$.
  Clearly $\Phi_\FG|_\Lambda$ is injective since $\Lambda$ fibers on 
  $\ell_1\cap\Lambda$, each fiber being a leaf of $cF_\xi$, with $G$ strictly monotonic
  on each fiber and $F$ strictly monotonic on the base. Moreover, 
  $F(\Lambda)\subset(-\infty,a)$ by assumption.
  Similarly, each $A_i$ fibers on $\ell_i\cap A_i$ so that $\Phi|_{A_i}$ is
  injective too; this time though $F(A_i)\subset[a,\infty)$ and, moreover, $G(A_i)=G(s_i)$.
  Consider now the set $V'=V\cap G^{-1}((c_1,c_2))$ and let $s$ be any leaf 
  of $\cF_\xi$ inside $\Lambda$. The sets of leaves of $\cG|_{V'}$ intersecting, 
  respectively, $s_1$ and $s_2$ cut $s$ in two disjoint open intervals
  $(c_1,b_1)$ and $(b_2,c_2)$; in particular all leaves of $\cG|_{V'}$ corresponding
  to the values in the closed interval $[b_1,b_2]$ do not intersect neither
  $s_1$ nor $s_2$ and are such that $s_1$ and $s_2$ lie on different components 
  with respect to each of them.
  Finally, let $F(\ell_1)=(a_1,a'_2)$ and $F(\ell_2)=(a_1,a''_2)$.
  Then $\Phi_\FG(\Lambda\cap G^{-1}((c_1,c_2)))=(a_1,a)\times(c_1,c_2)$,
  $\Phi_\FG(A_1\cap G^{-1}((c_1,c_2)))=[a,a'_2)\times(c_1,b_1)$
  and
  $\Phi_\FG(A_2\cap G^{-1}((c_1,c_2)))=[a,a''_2)\times(b_2,c_2)$
  so that $\Phi_\FG(V\cap F^{-1}((a_1,a_2))\cap G^{-1}((c_1,c_2)))=(a_1,a_2)\times(c_1,c_2)\setminus R$
  for $a_2=\min\{a'_2,a''_2\}$. 

  In case $V$ is contained in $F>a$, we use the chart $\tilde \Phi_\FG=(-F,G)$ and 
  repeat the argument above.
\end{proof}
%
We call the chart $(U\cap \Phi_\FG^{-1}((a_1,a_2)\times(c_1,c_2)),\Phi_\FG)$ 
\footnote{Replace $\Phi_\FG$ with $\tilde\Phi_\FG$ if, in the terminology of Proposition~\ref{thm:nc},
$V$ is contained in $F>a$.} a \emph{normal chart for $s_1$ and $s_2$}.
By Proposition~\ref{thm:smoothSols} there are countably many conditions that 
must be satisfied for each one of the intervals between pairs of adjacent
separatrices so that equation~(\ref{eq:Fy=G}) admits a smooth solution. 
Since in $\Gamma_\FG$ there is a natural family of transversals for $\cF_\xi$
these conditions can be restated more properly for this setting in the following way.
Let $I=\{a\}\times[b_1,b_2]$ the vertical interval separating a pair of adjacent 
separatrices $s_1$ and $s_2$ in a normal chart.
Every such interval determines a rings homomorphism
$\theta^{(r)}_I:SG^r_I\to SG^r_a$
defined as follows. Given $\fg\in SG^r_I$, let $\hat g\in\fg$ and $\delta=\min\{c_2-b_2,b_1-c_1\}$, 
choose an arbitrary $\epsilon\in(0,\delta)$ and set 
$h_I(x') = \int_{b_1-\epsilon}^{b_2+\epsilon}\hat g(x',y')dy'$ for $x'\in(a_1,a)$; we
define $\theta^{(r)}_I(\fg)=[h_I]_{SG_a^r}$.
\begin{proposition}
  The left singular germ of $h_I$, modulo germs of smooth functions, does not depend on 
  the particular choice of $\epsilon\in(0,\delta)$ and $\hat g\in\fg$.
\end{proposition}
\begin{proof}
  Let $h'_I(x') = \int_{b_1-\epsilon'}^{b_2+\epsilon'}\hat g'(x',y')dy'$ for 
  $\epsilon'\in(0,\delta)$ and $\hat g'\in\fg$. Then the function
  $$
  h'_I(x') - h_I(x') = 
  \int_{b_1-\epsilon'}^{b_2-\epsilon}\left(\hat g'(x',y')-\hat g(x',y')\right)dy'+
  $$
  $$
  +
  \int_{b_2+\epsilon}^{b_2+\epsilon'}\hat g(x',y')dy'
  +
  \int_{b_1-\epsilon'}^{b_1-\epsilon'}\hat g(x',y')dy'
  $$
  is smooth in $(a_1,a]$ since the integrands are all smooth in 
  $R_I$, the last two because $g$ is smooth in $R_I\setminus I$ and the integral
  intervals lie inside that set for every $x\in(a_1,a]$ and the first 
  because by hypothesis $\hat g'-\hat g$ is identically zero in some left 
  neighbourhood of $I$. Adding to $g$ and $\hat g$ any function smooth in 
  the whole $R_I$ changes the rhs just by a smooth function. 
\end{proof}
The maps $\theta^{(r)}_I$ then are well-defined.
It is clear from the definition of $h_I$ that
$\theta^{(r)}_I$ is a $C^r_x(\bR)$-module homomorphism, where $C^r_x(\bR)$ is the algebra
of $C^r$ functions depending on $x'$ only, since 
$$
\int_{b_1-\epsilon}^{b_2+\epsilon}f(x')\hat g(x',y')dy'=f(x')\int_{b_1-\epsilon}^{b_2+\epsilon}\hat g(x',y')dy'\,,
$$
and commutes with the derivatives with respect to $x'$, i.e.
$\theta^{(r)}_I(\partial^k_{x'}\hat g)=\partial^k_{x'}\theta^{(r)}_I(\hat g)$.

Next proposition shows that the maps $\theta^{(r)}_I$ determine the solvability of
the cohomological equation.
%
\begin{thmGerms}
  Let $\{I_j\}$ be the set of all (vertical, closed) intervals between adjacent 
  separatrices in $\Gamma_\FG$ and $\theta^{(r)}_{I_j}$ the corresponding ring homomorphisms.
  Then $g\in L_{\xi'_\F}(C^r(\Rt))$ iff $[(\hat\Phi_\FG)_* g]_{SG^r_{I_j}}\in\ker\theta^{(r)}_{I_j}$
  for all $\theta^{(r)}_{I_j}$.
\end{thmGerms}
\begin{proof}
  Let $I=\{a\}\times[b_1,b_2]$ be the vertical interval which separates two 
  adjacent separatrices of $\xi'_\F$ in a normal chart for the corresponding
  adjacent separatrices $s_1$ and $s_2$ 
  and set $\hat g=(\Phi_\FG)_* g$ within that chart. 
  Then 
  $$
  \lim_{x'\to a^-}\int_{b_1-\epsilon}^{b_2+\epsilon}\partial^k_{x'}\hat g(x',y')dy'
  $$ 
  is exactly the gap of $\Phi_\FG^*g$ between $s_1$ and $s_2$ with respect to the pair 
  of transversals which are the counterimages of $y'=b_1-\epsilon$ and $y'=b_2+\epsilon$ 
  and the gap exists and is finite if and only if those functions can all be extended 
  to continuous
  functions for all $k<r$, which in turn means that the (germ of the) function
  $\int_{b_1-\epsilon}^{b_2+\epsilon}\partial^k_{x'}\hat g(x',y')dy'$ can be extended to 
  a smooth map up to $x'=a$, i.e. $[(\hat\Phi_\FG)_* g]_{SG^r_I}\in\ker\theta^{(r)}_I$.
  Now the claim follows immediately from Theorem~\ref{thm:smoothSols}.
\end{proof}
%
The $C^r(\bR)$-modules $\Theta^r_{I_j}=\ker\theta^{(r)}_{I_j}$ contain therefore the (left singular) 
germs of all functions for which the cohomological equation is solvable in the 
neighbourhood of a pair of adjacent separatrices.
Modulo isomorphisms there are only two such spaces:
the one relative to $J=\{0\}\times[-1,1]$ and the one relative to $O=\{(0,0)\}$;
moreover $\Theta^r_O\subset\Theta^r_J$.
%
\begin{proposition}
  \label{thm:Theta}
  The spaces $\Theta^r_O$ satisfy the following properties:
  \begin{enumerate}
    \item $\Theta^r_O$ contains the singular left germs of all $y'$-odd\footnote{We say 
      that $f(x,y)$ is $y$-odd if $f(x,-y)=-f(x,y)$ and $y$-even if $f(x,-y)=f(x,y)$.} 
      $C^r$ functions;
    \item $\Theta^r_O$ contains the singular left germs of some but not all $y'$-even 
      $C^r$ functions;
    \item $\Theta^{r+1}_O$ is strictly contained in $\Theta^r_O$.
  \end{enumerate}
\end{proposition}
\begin{proof}
  1. If $\hat g$ is $y'$-odd then also every $\partial^k_{x'}\hat g$ is so for every 
  $k\leq r$; then $\int_{-\epsilon}^{\epsilon}\partial^k_{x'}\hat g(x',y')dy'$ 
  is identically zero for every $k\leq r$ and therefore it can be extended smoothly 
  to a $C^r$ function up to $x'=0$.

  2. Consider $\hat g(x',y')=e^{-(y')^2/(x')^2}/\sqrt{-\pi x'}\in C^\infty(\Rt\setminus(0,0))$, 
  so that
  $$
  \lim_{x\to0^-}g(x',y')=0\,,\;y'\neq0\,;\;\lim_{x\to0^-}g(x',0)
  =\infty\,;\;\int_{-\infty}^{\infty}g(x',y')dy'=1\,,\;\forall x'\in\bR\,.
  $$
  By reparametrizing the $y'$ coordinate we can find a $\hat g'$ with the same limits with 
  respect to $x'\to0$ but such that $\int_{-\epsilon}^{\epsilon}\hat g'(x',y')dy'=1$. 
  Since the $\theta^{(r)}_O$ are homomorphisms of $C^r_x(\bR)$-modules we can get in this 
  way every $C^r$ function $f(x')$ just by multiplying $g'(x',y')$ by $f(x')$. 
  On the other side, germs of functions diverging too fast, e.g. as 
  $\hat g(x',y')=(x')^{-2}+(y')^{-2}$, do not belong to any $\Theta^r_O$.

  3. Consider $\displaystyle\hat g(x',y')=\frac{x'}{\sqrt{(x')^2+(y')^2}}\in C^\infty(\Rt\setminus(0,0))$.
  The germ of the corresponding 
  $h_O(x')=2x'\log\left[2\left(y'+\sqrt{(x')^2+(y')^2}\right)\right]_{y'=0}^{\epsilon}$
  can be extended at 0 to a $C^0$ (but not $C^1$) function. By integrating $r$ times $\hat g$ with
  respect to $x'$ one can get concrete examples of functions smooth in $\Rt\setminus(0,0)$
  whose germ belongs to $\Theta^r_O$ but not to $\Theta^{r+1}_O$.
\end{proof}
An immediate consequence of point (3) of the proposition above is the following:
\begin{corollary}
  Let $\xi\in\vf_r(\Rt)$, $L^{(r)}_\xi$ the restriction of $L_\xi$ to $C^r(\Rt)$ and let 
  $C^r_\xi(\Rt)$ be the set of all functions $f\in C^r(\Rt)$
  such that $f+g$ is at most $C^r$ for all $g\in\ker L^{(r)}_\xi$. The inclusions 
  $$
  L^{(r+1)}_\xi\left(C^{r+1}_\xi(\Rt)\right)\cap\smooth\subset L_\xi^{(r)}\left(C^{r}_\xi(\Rt)\right)\cap\smooth
  $$ 
  are strict for every $r\in\bN$.
\end{corollary}
\begin{proof}
  The fact that 
  $L^{(r+1)}_\xi\left(C^{r+1}(\Rt)\right)\cap\smooth\subset L_\xi^{(r)}\left(C^{r}(\Rt)\right)\cap\smooth$
  is trivially true because $L_\xi^{(r)}(f+g)\in\smooth$ for each $f\in\smooth$, $g\in\ker L^{(r)}_\xi$.
  Our claim is that the inclusion is true even when we restrict $L_\xi$ to the space of functions which
  are ``strongly $C^r$'' with respect to $\xi$, i.e. those that cannot be made smoother by adding
  to them an element of the kernel of $L^{(r)}_\xi$. Consider indeed the concrete case used in point (3)
  of Proposition~\ref{thm:Theta}: in a normal chart, where the two separatrices are given by
  $x'=0$, $y'>a$ and $x'=0$, $y'<a$, the (local) primitive of $\hat g(x',y')=x'/\sqrt{(x')^2+(y')^2}$ 
  is $f(x',y')=x'\log\left[2\left(y'+\sqrt{(x')^2+(y')^2}\right)\right]$, which is $C^0$ but not $C^1$
  because the first derivative with respect to $x'$ diverges on the second separatrix. Since the divergence
  takes place only on one of the separatrices, there is no way to eliminate it by adding a function belonging
  to the kernel of $L_\xi$.
\end{proof}
In the following subsections we work out in detail two model examples.
%
\subsection{$\xi_n=(1-n+(1+n)y)\,\partial_x+(1-y^2)\,\partial_y$}
%
The $\xi_n$, $n\in\bN$, are all of finite type since they are polynomial.
In particular they all have exactly two separatrices, the straight lines $y=\pm1$,
which bound the canonical region $\bR\times(-1,1)$. 
The function $F_n(x,y)=(1+y)^n(1-y)e^x$ is a functional generator for $\ker L_{\xi_n}$
so the only intrinsically Hamiltonian among them is $\xi_1=2y\,\partial_x+(1-y^2)\,\partial_y$.
All of them are transversal to the same Hamiltonian foliation $\cG$ of the level sets 
of $G(x,y)=y e^x$, which is topologically conjugate with the trivial foliation in parallel 
straight lines. The 2-form 
$$
\Omega_\FG=2(1+y)^{n-1}(1 - (n-1)y + n y^2)e^{2x}\,\Omega_0
$$
is degenerate on the separatrix $y=-1$ except in the $n=1$ case, when is globally non-degenerate.
Via $\Omega_\FG$ we get 
$$
\xi'_{\Fn}=\frac{1}{2e^x(1 - (n-1)y + n y^2)}\xi_n\,,\;\xi'_{\Gn}=\frac{1}{2e^x(1+y)^{n-1}(1 - (n-1)y + n y^2)}\eta\,,
$$ 
where $\eta=2\partial_x-2y\partial_y$. Due to the degeneracy of $\Omega_\FG$, $\xi'_{\Gn}$ diverges
on the separatrix $y=-1$ for $n\neq1$.

The image of every $\Phi_\FnG$ is $\Rt_0=\Rt\setminus\{0\}\times[0,\infty)$ and $\Phi_\FoG$ 
is an almost complex map between $(\Rt,J_\FG)$ and $(\Rt_0,i)$ 
for
$$
J_\FG=y\,\partial_x\otimes dx + 2\,\partial_x\otimes dy - (1+y^2)/2\,\partial_y\otimes dx - y\,\partial_y\otimes dy\,.
$$
The leaves of $\cF_{\Fn}$ within
the canonical region are sent to the vertical lines of the half plane $x<0$ and the separatrices
$y=-1$ and $y=+1$ to the half lines $\{0\}\times(-\infty,0)$ and $\{0\}\times(0,+\infty)$ respectively.
The leaves lying in the half-plane $y>1$ fill in the vertical half-lines the first quadrant and 
the ones lying in $y<1$ the fourth quadrant. In this case the maps $\Phi_\FnG$ are all globally injective.
The cohomological equation $L_{\xi_\Fn}f=g$ maps to 
\begin{equation}
  \label{eq:p2}
  \partial_{y'}\hat f(x',y')=\hat g(x',y')\,,\;\hat g\in C^\infty(\Rt_0)\,.
\end{equation}
When $\hat g$ is smooth on the whole plane clearly (\ref{eq:p2}) is always solvable. E.g. all 
smooth solutions to 
$$
L_{\xi'_\Fn} f(x,y)=F_n(x,y) G(x,y)=2(y^2-1)(y+1)^{n-1}y e^{2x}
$$ 
are given by 
$$
f(x,y)=\frac{F_n(x,y) G^2(x,y)}{2} + h\left(F_n(x,y)\right) = 2(y^2-1)(y+1)^{n-1}y^2e^{3x}+h\left(F_n(x,y)\right)\,,
$$ 
where $h\in C^\infty(\bR)$.

In the following we assume $n=1$ since expressions are much simpler in this case.
Consider first the $y'$-odd function 
$$
\hat g(x',y')=\frac{y'}{\sqrt{(x')^2+(y')^2}}\in C^\infty(\bR_0)\,,\;\Phi^*_\FG \hat g(x,y)=\frac{2y}{1+y^2}\in\smooth\,.
$$
By Proposition~\ref{thm:Theta} the singular left germ of $\hat g$ belongs to 
$\Theta^\infty_O$ and therefore $g\in L_\xi(\smooth)$. 
Indeed (\ref{eq:p2}) in this case is solved by 
$$
\hat f(x',y')=\sqrt{(x')^2+(y')^2}\,,
$$ 
whose pull-back 
$$
\Phi^*_\FG f(x,y)=(1+y^2)e^x
$$ 
is globally smooth. Similarly, $y\in L_{\xi'_\F}(\smooth)$ since $y=\Phi_\FG^*\hat g(x,y)$
for the $y'$-odd singular function $\hat g(x',y')=(\sqrt{(x')^2+(y')^2}+x')/y'$.

On the contrary, in case of 
$$
\hat g(x',y')=\frac{x'}{\sqrt{(x')^2+(y')^2}}\,,\;g(x,y)=\Phi^*_\FG\hat g(x,y)=\frac{1-y^2}{1+y^2}\,,
$$ 
as discussed in Proposition~\ref{thm:Theta} we have that the germ of $\hat g$ belongs to 
$\Theta^0_O$ but not to $\Theta^1_O$; correspondingly all solutions will be $C^0$ but not $C^1$.
E.g. an explicit solution is given by
$$
f(x,y)=\Phi^*_\FG\left(x'\log\left[2\left(y+\sqrt{(x')^2+(y')^2}\right)\right]\right)=(1-y^2)e^x\left(x+2\log\left|1+y\right|\right)\,.
$$ 
Note that Lie derivatives of $f$ are, as expected, smooth with respect to $\xi'_\Fo$ direction 
but singular (on the horizontal straight line $y=-1$) with respect to $\eta$. 
In particular, $g$ belongs to $L_{\xi'_\F}(L^1_{\loc}(\Rt))$ (where the derivative is intended 
in the weak sense) but does not belong to any $L_{\xi'_\Fo}(C^k(\Rt))$, $k>1$.
The same happens in case of $x=\Phi_\FG^*\hat g(x,y)$, where 
$\hat g(x',y')=\log(\sqrt{(x')^2+(y')^2}+x')/2)$.
For a thorough discussion about locally integrable solutions of regular polynomial vector 
fields in the plane depending only on one variable see~\cite{DGK10}.
\subsection{$\xi_n=\left(\cos y+(n-1)\cos^2(y/2)\right)\partial_x-\sin y\,\partial_y$}
%
\begin{figure}
  \centering
  \includegraphics[width=5cm]{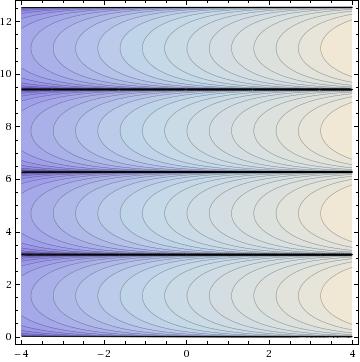}\hskip.75cm\includegraphics[width=5cm]{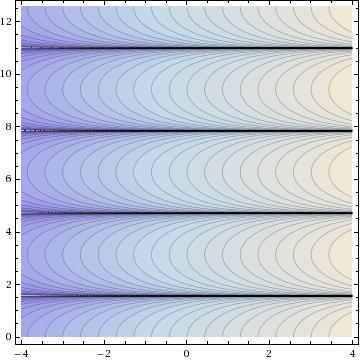}
  \caption{%
    \small
    Level sets of $F(x,y)=e^x\sin y$ (left) and $G(x,y)=e^x\cos y$ (right). The  
    separatrices of the first foliation are the straight lines $s_k=\{y=k\pi\}$, $k\in\bZ$, 
    the ones of the second are the straight lines $s'_k=\{y=\pi/2+k\pi\}$, $k\in\bZ$. 
    Note that $\cI_{s_n}=\{s_{n-1},s_{n+1}\}$, i.e. $s_n$ is inseparable only from $s_{n-1}$ 
    and $s_{n+1}$ (this is possible because the relation of inseparability is not 
    transitive). The same holds for the $s'_k$.
  }
  \label{fig:expz}
\end{figure}
The $\xi_n$, $n\in\bN$, are all of finite type for their components are Morse
functions depending only on one variable; in this case indeed only the vertical 
lines can be separatrices and they do not accumulate within any compact set.
For every $\xi_n$ the set of separatrices is $\cS=\{y=k\pi\,,\;k\in\bZ\}$. 
The function $F_n(x,y)=-\sin^{n-1}(y/2)\sin y\, e^x$ is a 
functional generator for $\ker L_{\xi_n}$ so that the only intrinsically Hamiltonian among 
them is $\xi_1=\cos y\,\partial_x-\sin y\,\partial_y$.
A Hamiltonian transversal foliation $\cG_n$ for $\xi_n$ is given by the level sets of
$G_n(x,y)=\cos y\, e^{x/n}$. The 2-form 
$$
\Omega_\FG=\left[(n-1)(2\cos y-\cos(2y))+3n+1\right]\sin^{n-1}(y/2)e^{(n+1)x/n}\,\Omega_0/4n
$$
is degenerate on the separatrices $y=2k\pi$, except of course in the $n=1$ case
when is globally non-degenerate. Via $\Omega_\FG$ we get 
$$
\xi'_{\Fn}=\frac{2n e^{-x/n}}{n+1+(n-1)(\sin^2y-\cos y)}\xi_n,\,
\xi'_{\Gn}=\frac{-2 \sin^{1-n}(y/2) e^{-x}}{n+1+(n-1)(\sin^2y-\cos y)}\eta\,,
$$ 
where $\eta=n\sin y\,\partial_x + \cos y\,\partial_y$. Due to the degeneracy of 
$\Omega_\FG$, $\xi'_{\Gn}$ diverges on the separatrices $y=2k\pi$, $k\in\bZ$, for $n\neq1$.

The image of every $\Phi_\FnG$ is $\Rt\setminus\{(0,0)\}$. Note that $\Phi_\FoGo$ is an 
almost complex map with respect to the almost complex structure
$$
J_\FoGo= \partial_y\otimes dx - \partial_x\otimes dy\,,
$$
so that $\Phi_\FoGo$ is actually a 
holomorphic map; in fact, in complex coordinates, $\Phi_\FoGo(z)=e^{z+i\pi/4}$ and its
graph is the Riemann surface of the complex logarithm. 
The graphs of all other $\Phi_\FnGn$ are diffeomorphic to it.

Consider just the case of the coordinate functions $x$ and $y$. 
The first is $y'$-even since $2x=\Phi_\FG^*\hat g(x,y)$ for 
$\hat g(x',y')=\log\left[(x')^2+(y')^2\right]$. A direct calculation shows that
$$
\left[\theta_n(\hat g)\right](x')=
2\int_0^\epsilon\log\left((x')^2+(y')^2\right)dy'=
4 x\tan^{-1}(\epsilon/x) + 2\epsilon (\log(\epsilon^2 + x^2)-2)
$$
which can be continued to a smooth function up to $x'=0$. 
Hence $\hat g\in\Theta^\infty_n$ for all $n$ and, correspondingly,
$x\in L_{\xi'_\F}(\smooth)$. An explicit solution is given by 
$$
f(x,y)=\Phi^*_\FG[2x'\tan^{-1}\frac{y'}{x'}+\log[(x')^2 + (y')^2]-2y']
=2\left[(x-1) \cos y - y \sin y \right]e^x\,.
$$
The second is $y'$-odd since $y=\Phi_\FG^*\hat g(x,y)$ for $\hat g(x',y')=\tan^{-1}(x'/y')$.
Hence even in this case $\hat g\in\Theta^\infty_n$ for all $n$, i.e. $y\in L_{\xi'_\F}(\smooth)$.
An explicit solution is given by 
$$
f(x,y)=\Phi^*_\FG[y' \tan^{-1}\frac{x'}{y'} + \frac{1}{2}\log[(x')^2 + (y')^2]]
=-\left[y \cos y + x \sin y\right]e^x\,.
$$
\section*{Acknowledgments}
The author wants to thank G.~D'Ambra and A.~Loi for introducing the subject, 
S.P.~Novikov for encouraging its study and for precious suggestions and 
A. Bergamasco, G.~D'Ambra, T. Gramchev, A.~Loi, R.~Mossa and D.~Zuddas 
for several fruitful discussions.
%
\bibliography{ceonr2}
\end{document}